\documentclass[journal ]{new-aiaa}
\usepackage[utf8]{inputenc}
\usepackage{textcomp}
\usepackage{mathtools}
\usepackage{graphicx}
\usepackage{amsmath}
\usepackage[version=4]{mhchem}
\usepackage{siunitx}
\usepackage{longtable,tabularx}
\usepackage[ruled,vlined]{algorithm2e}

\newcommand{\tr}{\text{{Tr}}}
\newcommand{\col}{\text{{Col}}}
\newcommand{\vect}{\text{{ vec}}}

\title{Control Allocation for Hybrid Coulomb Spacecraft Formations}

\author{Adam M Tahir\footnote{William E Boeing Department of Aeronautics and Astronautics, University of Washington, Seattle, WA.} }
\affil{University of Washington, Seattle, WA, 98115}

\begin{document}

\maketitle

\begin{abstract}
This paper proposes an algorithm which can be used in hybrid Coulomb spacecraft formations to minimize propellant by maximizing the amount of force that is generated by Coulomb forces. This problem is difficult due to the nonlinearities inherent in Coulomb's law. The problem is posed as a series of matrix rank minimization problems which can be solved efficiently using a trace heuristic. A numerical example is provided which shows that, using the proposed control allocation algorithm, the amount of propellant required to perform a reconfiguration maneuver is reduced by approximately 40\% compared to using solely thrusters. 
\end{abstract}


\section{Introduction}
The lifespan of a spacecraft formation is limited by the amount of propellant that can be brought onboard to be used for stationkeeping and other formation control maneuvers. As a measure to reduce propellant requirements for spacecraft formation control, Schuab {\it et al.} \cite{Prospects} first proposed the prospect of using the inter-spacecraft Coulomb forces for formation control in high orbits, which could be done with negligible amounts of propellant. Coulomb formation control also brings with it many technical challenges. The two main challenges are the nonlinearities and underactuation that inherent from Coulomb's law \cite{Prospects, CoulombNonlin}. Future missions that take advantage of Coulomb forces will likely have to be  hybrid Coulomb spacecraft formations to overcome these challenges. 

In a hybrid Coulomb spacecraft formation, each spacecraft is equipped with two methods of actuation:
\begin{enumerate}
\item {\bf Active Charge Control}. Each spacecraft can actively change its charge. By Coulomb's law, the charged spacecraft exert action-reaction pairs of forces on each other.  These forces are proportional to the product of their charges, inversely proportional to the square of the separation distance, and act along the line of sight between two spacecraft.
\item {\bf Thrusters.} Each spacecraft is equipped with a set of thrusters that can produce a thrust forces in any direction and any magnitude. 
\end{enumerate}
The thrusters complement the Coulomb actuation by producing forces that the Coulomb actuation cannot. The role of thrusters for a two-spacecraft hybrid Coulomb formation is straightforward: the thrusters are used to generate forces perpendicular to the line-of-sight between the two spacecraft (cf. \cite{ArunHybrid, CoulombNonlin}).  For hybrid formations of more than two spacecraft, the problem is more difficult and cannot be solved analytically except under strong assumptions such as symmetry.

In this paper, the control allocation problem for hybrid Coulomb spacecraft formations is defined and a control allocation algorithm is proposed that can be run efficiently. The efficiency of the proposed algorithm comes from its use of convex optimization. 










\subsection{Notation}
The $n\times n$ identity matrix is denoted by $I_n$. A matrix $P\in\mathbb{R}^{n\times n}$ is symmetric positive semidefinite if $x^\top Px\ge 0$ for all $x\in\mathbb{R}^{n}$.  The notation $P\succeq 0$ denotes that $P$ is symmetric positive semidefinite. For a matrix $A\in\mathbb{R}^{n\times n}$, the notation He$(A)=\frac{1}{2}(A+A^\top)$. The Euclidean norm of a vector $x\in\mathbb{R}^n$ is denoted by $\| x \|$.  The Frobenius norm of a matrix $A\in\mathbb{R}^{n\times m}$ is denoted $\|A\|_F$. Two vectors $x,y\in\mathbb{R}^n$ are said to be $\epsilon$-close if $\|x-y\|\le \epsilon$. The trace operator is denoted by $\tr$. The notation $\col(x_1,\dots, x_N)$ denotes the column vector composed of stacking $x_1$ on top of $x_2$, etc. The notation $\vect(A)$ denotes the stacking of the columns of matrix $A$ on top of each other starting from the left-most column. Coulomb's constant is $k_c=8.99\times 10^{9} N \text{m}^2/\text{C}^2$. The Kronecker product of $A$ with $B$ is denoted by $A\otimes B$. The Moore-Penrose psuedoinverse of matrix $B$ is denoted by $B^\dagger=B^\top(B B^\top)^{-1}$. 

\section{Hybrid Coulomb Spacecraft Formations}
\begin{figure}[!h]
\centerline{\includegraphics[width=0.3\columnwidth]{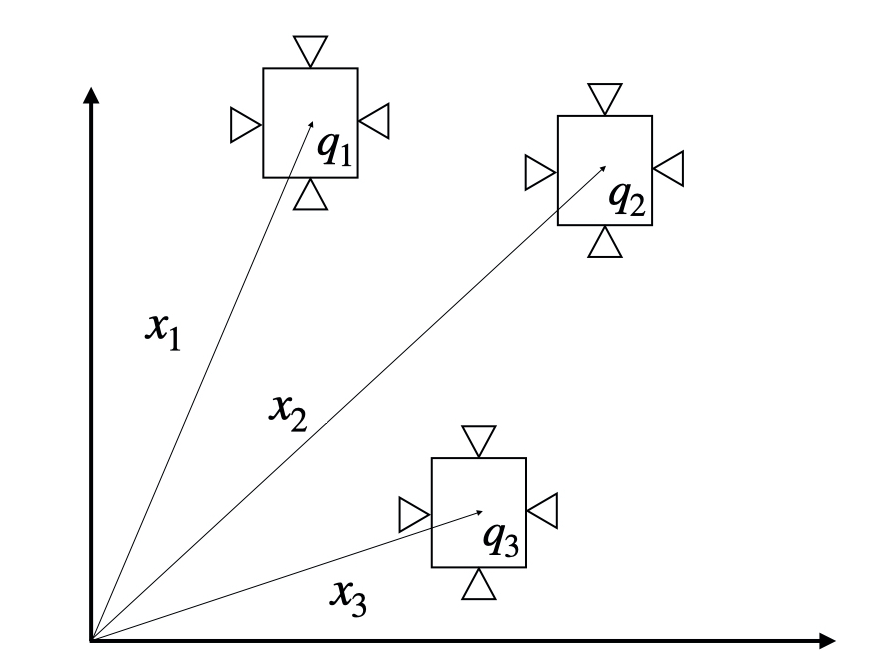}}
\caption{Schematic of a Coulomb spacecraft formation with $\mathcal{N}=3$ in $\mathbb{R}^2$ . The positions of each spacecraft are $x_1,x_2,x_3\in\mathbb{R}^2$, and their respective charges are $q_1,q_2,q_3\in\mathbb{R}$. Each spacecraft is depicted with a set of thrusters that can generate thrusts $T_1,T_2,T_3\in\mathbb{R}^2$.}
\label{f:CoulombSchematic}
\end{figure}

Consider a formation of $\mathcal{N}>1$ spacecraft embedded in $\mathbb{R}^d$. The position relative to an inertial frame of the spacecraft indexed $i\in\{1,\dots,\mathcal{N}\}$ is $x_i\in\mathbb{R}^d$, measured in m. Each spacecraft is modeled as a point charge with charge $q_i\in\mathbb{R}$, measured in C, which can be adjusted by its active charge control system. Moreover, each spacecraft is equipped with a set of conventional thrusters that can produce thrust $T_i\in\mathbb{R}^d$, measured in N. The vector of all positions in the formation is $x=\col(x_1,\dots, x_\mathcal{N})\in\mathbb{R}^{d\mathcal{N}}$, the vector of all charges is $q=\col(q_1,\dots, q_\mathcal{N})\in\mathbb{R}^\mathcal{N}$, and the vector of all thrusts is $T=\col(T_1,\dots,T_\mathcal{N})\in\mathbb{R}^{d\mathcal{N}}$. In Fig. \ref{f:CoulombSchematic}, a formation with $\mathcal{N}=3$ and $d=2$ is depicted. 

The force $F_{C_i}\in\mathbb{R}^d$, in N, acting on spacecraft $i$ is the vector sum of the inter-spacecraft Coulomb forces\footnote{Note that environmental effects such as Debye shielding are neglected.}:
\begin{align}
F_{C_i}(x,q) &=\sum_{j=1,j\ne i}^{\mathcal{N}} k_c\dfrac{x_i-x_j}{\|x_i-x_j\|^3}q_iq_j\nonumber \\
&= k_c a_i(x) vec(qq^\top),\label{e:Coulombforce}
\end{align}
where,
\begin{align}
a_i(x) = \sum_{j=1, j\ne i}^\mathcal{N}\left(\left(\vect\left(\text{He}(\Upsilon(i,j))\right)\right)^\top\otimes \left(\dfrac{x_i-x_j}{\|x_i-x_j\|^3}\right)\right),\label{e:aidef}
\end{align}
and $\Upsilon(i,j)\in\mathbb{R}^{\mathcal{N}\times\mathcal{N}}$ is a matrix which is zero everywhere except its $(i,j)$-th entry. Notice that the Coulomb forces depend on the charge of every spacecraft in the formation. 

 Since Coulomb forces are internal forces, the following holds:
$$\sum_{i=1}^\mathcal{N}F_{C_i}(x,q)=0,\; \forall (x,q)\in\mathbb{R}^{d\mathcal{N}}\times \mathbb{R}^\mathcal{N}.$$
Therefore, it is more convenient to work with {\it relative Coulomb forces}. The convention for relative Coulomb forces used in this paper is the following:
\begin{align}
\Delta{\bf F}_{C_i}(x,q) = F_{C_{i+1}}(x,q) -F_{C_i}(x,q), \;\forall i=1,\dots, \mathcal{N}-1.\label{e:Coulomb}
\end{align}
The vector of relative Coulomb forces is the following:
\begin{align}
\Delta{\bf F}_C(x,q) &=\col\left(\Delta{\bf F}_{C_1}(x,q) ,\dots,\Delta{\bf F}_{C_{\mathcal{N}-1}}(x,q) \right)\in\mathbb{R}^{d(\mathcal{N}-1)}\nonumber\\
&= k_c A(x)\vect(qq^\top),\label{e:DeltaFC}
\end{align}
where
\begin{align}
A(x) = \col\left(a_2(x)-a_1(x), \dots, a_{\mathcal{N}}(x)-a_{\mathcal{N}-1}(x)\right)\in\mathbb{R}^{d(\mathcal{N}-1)\times \mathcal{N}^2}\label{e:A}
\end{align}
and $a_i(x)$ is defined in \eqref{e:aidef}.

Using the relative force convention defined for Coulomb forces in \eqref{e:Coulomb}, the following is the vector of {\it relative thruster forces}:
\begin{align}
\Delta {\bf F}_T(T) &= \col(T_2-T_1,\dots, T_{\mathcal{N}}-T_{\mathcal{N}-1})\in\mathbb{R}^{d(\mathcal{N}-1)}\nonumber \\
&=BT,\label{e:relt}
\end{align}
where $B\in\mathbb{R}^{d(\mathcal{N}-1)\times d\mathcal{N}}$ is the following:
\begin{align}
B = \begin{bmatrix}
-1 & 1 & 0 & \dots & 0\\
0 & -1 & 1 &  \dots & 0\\
\vdots & & \ddots & \ddots &  & \\
0 &\dots & & -1 & 1
\end{bmatrix}\otimes I_d.\label{e:BB}
\end{align}
Note that $B$ has linearly independent rows and, thus, is right invertible (i.e. $BB^\dagger=I_{d(\mathcal{N}-1)}$). 
\section{Control Allocation Problem}

\begin{figure}[!ht]
\centerline{\includegraphics[width=0.8\columnwidth]{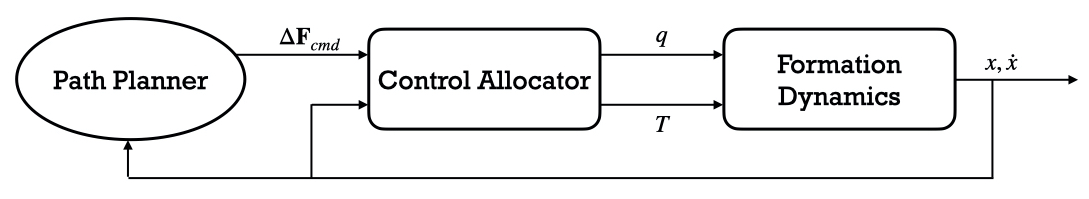}}
\caption{Control allocation problem.}
\label{f:CA}
\end{figure}

A control allocator receives a relative force command $\Delta {\bf F}_{\text{cmd}}=\col\left(\Delta {\bf F}_{{\text{cmd}}_1}, \dots, \Delta {\bf F}_{{\text{cmd}}_{\mathcal{N}-1}}\right)\in\mathbb{R}^{d(\mathcal{N}-1)}$ from a path planner and, based on the current formation geometry $x$, computes $q\in\mathbb{R}^\mathbb{N}$ and $T\in\mathbb{R}^{d\mathcal{N}}$ such that:
\begin{equation}
\Delta{\bf F}_T(T)+\Delta {\bf F}_C(x,q)=\Delta {\bf F}_{\text{cmd}},\label{e:CAF}
\end{equation}
where  $\Delta {\bf F}_C(x,q)$ and $\Delta{\bf F}_T(T)$ are defined in \eqref{e:DeltaFC} and \eqref{e:relt}, respectively.  This workflow is illustrated in Fig. \ref{f:CA}. 

The hybrid Coulomb spacecraft formation is over-actuated in the sense that any commanded relative force could be generated solely by thrusters  (i.e. \eqref{e:CAF} is always satisfied with $q=0$ and $T=B^\dagger\Delta {\bf F}_{\text{cmd}}$). That, however, is inefficient. Changing the charge of each of the spacecraft requires negligible amounts of propellant \cite{Prospects} and conventional thrusters require large amounts of propellant. The objective of the control allocator is, therefore, to generate as much of the relative force command from the Coulomb forces as possible to minimize the amount of propellant required to accomplish formation maneuvers. Ideally, this can be solved by the following optimization problem:
\begin{subequations}\label{e:genericopt}
\begin{align}
\underset{T\in\mathbb{R}^{d\mathcal{N}}, q\in\mathbb{R}^\mathcal{N}}{\text{minimize}} \; \;\;
& \|T\| \\
\text{subject to}\;\;\;
&  \eqref{e:CAF}
\end{align}
\end{subequations}
However, the optimization problem \eqref{e:genericopt} is difficult to solve directly since the equality constraint \eqref{e:CAF} is nonconvex  in the variable $q$.  In the sequel, an algorithm is proposed which uses a convex program to efficiently find approximate solutions to the control allocation problem. The convex program is posed so that a feasible solution can be used to find a set of charges $q$ that generates as much of the commanded relative force $\Delta {\bf F}_{\text{cmd}}$ through Coulomb actuation as possible and then computes $T$ to fulfill the remainder of the force command.  

\section{Proposed Control Allocation Algorithm}
\subsection{Algorithm Overview}
The proposed control allocator shown in Algorithm \ref{algo:bro} solves the following convex semidefinite program for values of $\epsilon$ from a user-defined finite search set $\mathcal{E}\subset\left[0,\|\Delta {\bf F}_{\text{cmd}}\|\right)$: 
\begin{subequations}\label{e:trace}
\begin{align}
\underset{\mathcal{Q}\in\mathbb{R}^{\mathcal{N}\times \mathcal{N}}}{\text{minimize}} \; \;\;
& \tr (\mathcal{Q})\\
\text{subject to}\;\;\;
&  \left\|A(x)\vect(Q)-\Delta{\bf F}_{\text{cmd}}\right\|\le \epsilon \\
& \mathcal{Q}\succeq 0
\end{align}
\end{subequations}
where $A(x)$ is defined in \eqref{e:A}.  For each value of $\epsilon$, a viable set of charges $\tilde{q}$ is extracted from an eigenvalue decomposition of the $\mathcal{Q}$ that optimizes \eqref{e:trace}. From this viable $\tilde{q}$, a set of thrusts $\tilde{T}$ is computed to complete the commanded relative force. The output of the algorithm is the $\tilde{T}$ and $\tilde{q}$ combination that has the minimum $\|\tilde{T}\|$. 

\begin{algorithm}[H]\label{algo:bro}
\SetAlgoLined
\KwData{$x, \Delta {\bf F}_{\text{cmd}}, \mathcal{E}$}
 \Begin{
Initialize $T\to B^\dagger \Delta {\bf F}_{\text{cmd}}$ and  $q \to 0$\\
 
 \For{$\epsilon\in\mathcal{E}$}    
        { 
        		Solve \eqref{e:trace} for $\mathcal{Q}$\\
 	Compute the largest eigenvalue of $\mathcal{Q}$, $\lambda_\mathcal{N}(Q)$, and its associated eigenvalue $v_\mathcal{N}\in\mathbb{R}^{\mathcal{N}}$\\	
	 $\tilde{q}\to\sqrt{\frac{\lambda_\mathcal{N}(Q)}{k_c}} v_\mathcal{N}$\\
 $\tilde{T}\to  B^\dagger\left(\Delta{\bf F}_{\text{cmd}}-\Delta {\bf F}_C(x,\tilde{q})\right)$ \\
	
	\If{$\|\tilde{T}\|\le \|T\|$}
	{
	$q\to \tilde{q}$ \\
	$T\to \tilde{T}$
	}
   }
\KwRet $q,T$ 
 }
 \caption{Control Allocator for Hybrid Coulomb Spacecraft Formations}
\end{algorithm}

The remainder of this section is dedicated to explaining the proposed control allocator in Algorithm \ref{algo:bro} and the associated convex optimization problem \eqref{e:trace}.

\subsection{Analysis of Algorithm \ref{algo:bro}}
\subsubsection{The Optimization Problem \eqref{e:trace}}
The proposed control allocator is based on the following variable substitution\footnote{Coulomb's constant $k_c$ is included as a scaling factor so that the values of $q$ that are extracted using the variable substitution are on the order of 1-100's of $\mu$C.}:
\begin{align}
\mathcal{Q}= k_c qq^\top\in\mathbb{R}^{\mathcal{N}\times\mathcal{N}}, \label{e:varsub}
\end{align}
which satisfies the following\footnote{Disregarding the special case where rank$(\mathcal{Q})=0$ when $q=0$.}:
\begin{subequations}\label{e:constraintsQ}
\begin{align}
&\text{rank}(\mathcal{Q})=1,\label{e:rankconst}\\
&\mathcal{Q}\succeq 0.
\end{align}
\end{subequations}
The inverse of the variable substitution \eqref{e:varsub} is given as follows:
\begin{align}
q = \sqrt{\frac{\lambda_\mathcal{N}(Q)}{k_c}} v_\mathcal{N},\label{e:invQ}
\end{align}
where $\lambda_\mathcal{N}(\mathcal{Q})$ is the largest eigenvalue of $\mathcal{Q}$ and  $v_\mathcal{N}$ is its associated eigenvector. Note that \eqref{e:constraintsQ} implies that $v_\mathcal{N}$ is a real vector, and $\lambda_\mathcal{N}(\mathcal{Q})$ is real and positive.

The variable substitution \eqref{e:varsub}  is defined so that the relative Coulomb force can be written as linearly dependent on $\mathcal{Q}$ as follows:
\begin{align}
\Delta {\bf F}_C(x,\mathcal{Q}) = A(x)\vect(\mathcal{Q}),\label{e:substCoulomb}
\end{align}
but, to take advantage of the linearity of \eqref{e:substCoulomb} in $\mathcal{Q}$ in an algorithm,  the conditions \eqref{e:constraintsQ} must also be imposed. That is, suppose an algorithm, given a formation geometry $x$ and  commanded relative force $\Delta{\bf F}_{\text{cmd}}$, computes a matrix $\mathcal{Q}$ that satisfies \eqref{e:constraintsQ} and is such that the following holds for some $\epsilon\ge 0$:
\begin{align}
\left\|A(x)\vect(Q)-\Delta{\bf F}_{\text{cmd}}\right\|\le \epsilon. \label{e:epsclose}
\end{align}
The relative Coulomb force generated by the set of charges $q$ computed by \eqref{e:invQ} is guaranteed to be $\epsilon$-close to the $\Delta{\bf F}_{\text{cmd}}$. In contrast, if an algorithm computes a matrix $\mathcal{Q}$ that satisfies \eqref{e:epsclose} for some $\epsilon\ge 0$, but does not satisfy \eqref{e:constraintsQ}, the $\epsilon$-closeness of the Coulomb force with $q$ computed by \eqref{e:invQ} to the commanded force is not guaranteed. 

Although the $\epsilon$-closeness is not guaranteed when $\mathcal{Q}$ is not rank 1, if $\mathcal{Q}\succeq 0$, then $\tilde{\mathcal{Q}} = \lambda_\mathcal{N}(\mathcal{Q}) v_\mathcal{N}v_\mathcal{N}^\top$ is a best rank-one approximation to $\mathcal{Q}$ in the sense that it is a rank-one matrix that minimizes $\|\tilde{\mathcal{Q}}-\mathcal{Q}\|_F$ (cf. \cite[\S 7.4.2]{HornyJohnson}). Therefore, choosing $q$ according to \eqref{e:invQ} corresponds to finding a best rank-one approximation to $\mathcal{Q}$ and then computing the $q$ from the eigenvalue decomposition of that approximation. It is more likely that choosing $q$ according to \eqref{e:invQ} yields a relative Coulomb force that is $\epsilon$-close to $\Delta {\bf F}_{\text{cmd}}$ with a lower rank (i.e. ``closer to'' rank-one) symmetric positive semidefinite matrix $\mathcal{Q}$ that satisfies \eqref{e:epsclose}. Thus, a strategy emerges for a control allocator: try to minimize the rank of $\mathcal{Q}\succeq 0$ subject to \eqref{e:epsclose}, then find charges $q$ via \eqref{e:invQ}. Consider the following optimization problem as an {\it ansatz}:
\begin{subequations}\label{e:nonconvex}
\begin{align}
\underset{\mathcal{Q}\in\mathbb{R}^{\mathcal{N}\times \mathcal{N}}}{\text{minimize}} \; \;\;
& \text{rank}(\mathcal{Q})\\
\text{subject to}\;\;\;
&  \left\|A(x)\vect(\mathcal{Q})-\Delta{\bf F}_{\text{cmd}}\right\|\le \epsilon\\
& \mathcal{Q}\succeq 0
\end{align}
\end{subequations}
 The rank of a matrix is nonconvex, so \eqref{e:nonconvex} cannot be solved efficiently. A widely used heuristic to minimize rank of a matrix is to minimize the nuclear norm of a matrix, which is a convex function  (see \cite{GoingNuclear,NoiseyMatrixCOmpletion} and others). For symmetric positive semidefinite matrices, minimizing the nuclear norm is equivalent  to minimizing the trace of the matrix.  Therefore, the optimization problem \eqref{e:trace} with the trace objective  that is used in Algorithm \ref{algo:bro} is a convex relaxation of the nonconvex problem \eqref{e:nonconvex} that promotes low rank $\mathcal{Q}$ and can be solved efficiently multiple times for different values of $\epsilon$. Solving \eqref{e:trace} for different values of $\epsilon$ will find several viable sets of charges $q$. The criterion used to choose between the various viable sets of charges is the amount of thrust required.  
 \subsubsection{Computing Thrusts for a Set of Charges}
 Once a viable $q$ is found from $\mathcal{Q}$ from \eqref{e:trace}, the relative thrust $\Delta {\bf F}_T(T)$ is found using \eqref{e:CAF}. Then the thrusts of the individual spacecraft are computed by
\begin{align}
T=B^\dagger\left(\Delta{\bf F}_{\text{cmd}}-\Delta {\bf F}_C(x,q)\right),\label{optimalthrust}
\end{align}
where $B$ is defined in \eqref{e:BB}. The thrusts $T$ are chosen according to \eqref{optimalthrust} since that is the minimizer of $\|T\|^2$ subject to  $\Delta {\bf F}_T(T)=BT=\Delta{\bf F}_{\text{cmd}}-\Delta {\bf F}_C(x,q)$  \cite[\S 2.1]{Allocation}.

\subsubsection{The Search Set $\mathcal{E}$}\label{s:Search}
It is clear that if $\epsilon\ge \|\Delta {\bf F}_{\text{cmd}}\|$, then the optimal solution to \eqref{e:trace} is $\mathcal{Q}=0$. Hence, the search set should be restricted to the interval $\left[0,\|\Delta {\bf F}_{\text{cmd}}\|\right)$. 

\section{Numerical Examples}
\subsection{Analysis of Algorithm \ref{algo:bro}: Four Spacecraft Formation in $\mathbb{R}^2$}
Consider four spacecraft with the following positions:
\begin{align*}
x_1 = \col(0, 0), \; x_2 = \col(10, 0),\;x_3 = \col(5,7),\; x_4 = \col(-10,2),
\end{align*}
and consider the following relative force command:
\begin{align*}
\Delta {\bf F}_{\text{cmd}} = \col(-0.023,-0.067,-0.069,-0.211,-0.037, 0.1806).
\end{align*}

Using Algorithm \ref{algo:bro}, the following spacecraft charges are found: $q_1 = 36.61\mu$C, $q_2 = 19.56\mu$C, $q_3 = -27.08\mu$C, $q_4 = 16.25\mu$C. It can be clearly seen in Table \ref{table} that the Algorithm \ref{algo:bro} leads to significantly reduced thrust requirement compared to using solely thrusters (i.e. letting $q=0$ and $T=B^\dagger\Delta {\bf F}_{\text{cmd}}$). In total, Algorithm \ref{algo:bro}  reduces the required propellant by 82\% compared to solely using thrusters.  
 \begin{table}[!h]
\caption{Thrust required in Newtons to achieve $\Delta {\bf F}_{\text{cmd}}$ with solely thrusters and using Algorithm \ref{algo:bro}.}
\label{table}
\centering
\begin{tabular}{|c|l|l|}
\hline
$i$ & $T_i$ (Only Thrusters) & $T_i$ (Algorithm \ref{algo:bro}) \\
 \hline
 1 & $\col(0.0610, 0.1106)$ & $\col( -0.0049,-0.0227)$\\
 \hline
 2 & $\col(0.0380,0.0436)$ & $\col(  -0.0040, 0.0081)$\\
 \hline
  3 & $\col( -0.0310, -0.1674)$ & $\col(  -0.0166,   0.0120)$\\
 \hline
  4 & $\col( -0.0680, 0.0132)$ & $\col(  0.0255,0.0026)$\\
 \hline
\end{tabular}
\label{tab1}
\end{table}

\begin{figure}[!h]
\centerline{\includegraphics[width=0.5\columnwidth]{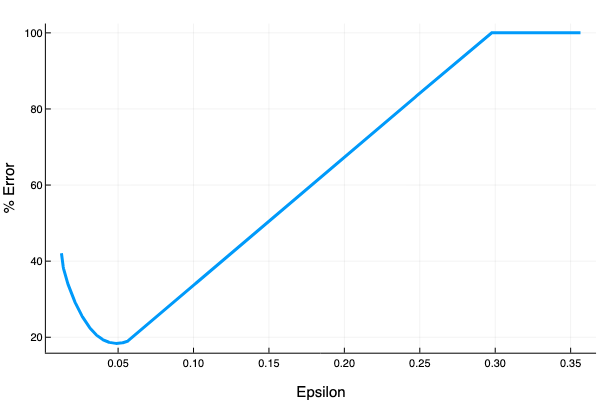}}
\caption{Percentage error of fit for different values of $\epsilon$.}
\label{f:fit}
\end{figure}

As part of Algorithm \ref{algo:bro}, the semidefinite program \eqref{e:trace} is solved for different values of $\epsilon$ using \verb|Convex.jl| \cite{JuliaConvex}. Fig. \ref{f:fit} plots the percentage error defined below:
\begin{align}
\text{\% Error} =\dfrac{\|\Delta{\bf F}_C(x,\tilde{q})-\Delta {\bf F}_{\text{cmd}}\|}{\|\Delta {\bf F}_{\text{cmd}}\|}\times 100\%\label{e:perc}
\end{align}
for the charge $\tilde{q}$ computed for different values of $\epsilon$. The iteration with $\epsilon=0.05$ produces charges with Coulomb force closest to the commanded relative force. 

\begin{figure}[!h]
\centerline{\includegraphics[width=0.5\columnwidth]{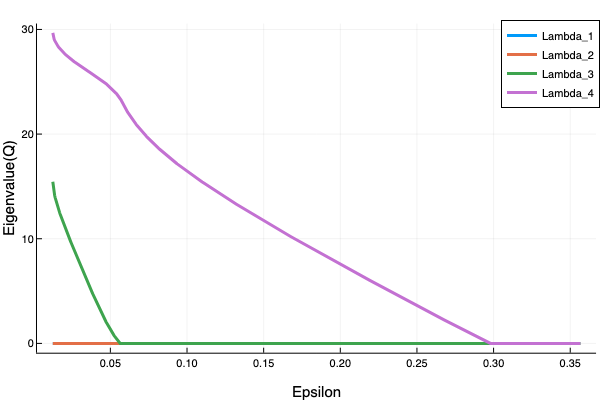}}
\caption{Eigenvalues of $\mathcal{Q}$ for different values of $\epsilon$.}
\label{f:eigens}
\end{figure}

Fig. \ref{f:eigens} plots the eigenvalues $\lambda_1,\dots,\lambda_4$ of the optimal $\mathcal{Q}$ for different values of $\epsilon$. 
The following observations can be made from Fig. \ref{f:eigens}:
\begin{enumerate}
\item  The optimal $\mathcal{Q}$ is a rank 1 matrix for all $0.055\le \epsilon\le 0.2971$, which can be seen from the fact that only one eigenvalue is nonzero for all  $0.055\le \epsilon\le 0.2971$. This shows how well the trace heuristic works for minimizing the rank. 
\item  As $\epsilon$ decreases below $\epsilon=0.05$, the second largest eigenvalue of $\mathcal{Q}$ increases. Meaning that the matrix $\mathcal{Q}$ is moving ``further away'' from the set of rank-one matrices. This contributes to the increase in the error seen in Fig. \ref{f:fit}  in this regime. This suggests that $\epsilon=0.05$ produces the best possible set of charges.
\item When $\epsilon\ge \|\Delta {\bf F}_{\text{cmd}}\|=0.2971$, the optimal $\mathcal{Q}$ is the zero matrix. This makes sense given the discussion in \S \ref{s:Search}. 
\end{enumerate}

\subsection{Propellant Savings for a Three-Spacecraft Reconfiguration Manuever in $\mathbb{R}^3$}
A reconfiguration maneuver is implemented for a formation with $\mathcal{N}=3$ in $\mathbb{R}^3$ with deep space dynamics. Algorithm \ref{algo:bro} is used to determine the charges and thrusts necessary to complete the maneuver while minimizing propellant required. 

Similarly to the relative force convention defined for Coulomb forces in \eqref{e:Coulomb}, the relative spacecraft positions are defined as 
\begin{align*}
\xi_i =x_{i+1}-x_i,\; \forall i=1,2.
\end{align*}

The reconfiguration maneuver is designed to drive the spacecraft from their initial relative configuration to a desired relative configuration: $\xi_i\to \xi_{\text{des}_i},\;  \forall i=1,2,$ where $\xi_{\text{des}_1}=\col(5,50,75),\; \xi_{\text{des}_2}=\col(60,25,100)$ and the spacecraft start from $\xi_1(0) =\col(100,0,0),\; \xi_{2}(0)=\col(0,0,100). $ The relative force commands are as follows:
\begin{align}
\Delta {\bf F}_{{\text{cmd}}_i} = -m\kappa(\xi_i-\xi_{\text{des}_i})-m\rho\dot{\xi}_i,\label{e:command_hist}
\end{align}
where $m=1$kg is the mass of all three spacecraft, and $\kappa=0.05,\varrho=0.2$ are chosen such that 
\begin{align*}
\begin{bmatrix}0 & I_3\\ -m\kappa I_3 & -m\varrho I_3 \end{bmatrix}\text{is Hurwitz.}
\end{align*}

\begin{figure}[!h]
\centerline{\includegraphics[width=0.5\columnwidth]{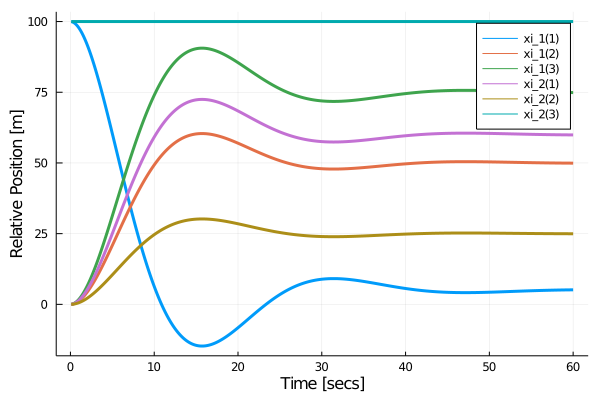}}
\caption{Trajectory of relative positions during the reconfiguration manuever.}
\label{f:traj}
\end{figure}

\begin{figure}[!h]
\centerline{\includegraphics[width=0.5\columnwidth]{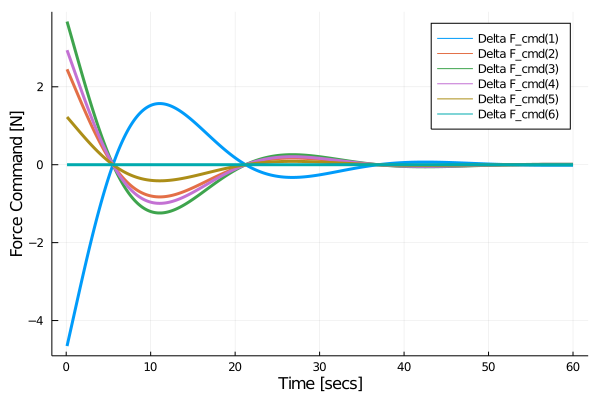}}
\caption{$\Delta {\bf F}_{\text{cmd}}$ history during the reconfiguration maneuver.}
\label{f:f_comtraj}
\end{figure}

\begin{figure}[!h]
\centerline{\includegraphics[width=0.5\columnwidth]{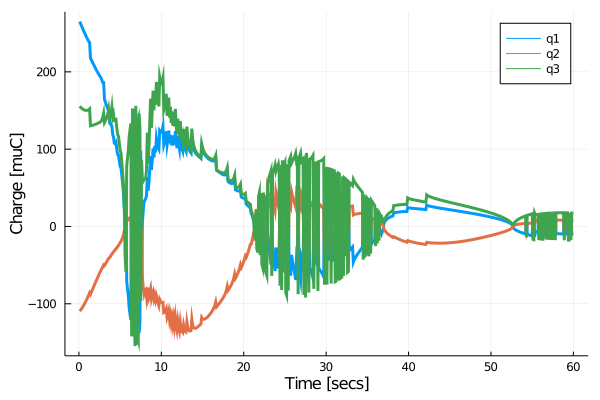}}
\caption{$q$ history during the reconfiguration maneuver computed by Algorithm \ref{algo:bro} .}
\label{f:charge_traj}
\end{figure}

\begin{figure}[!h]
\centerline{\includegraphics[width=0.5\columnwidth]{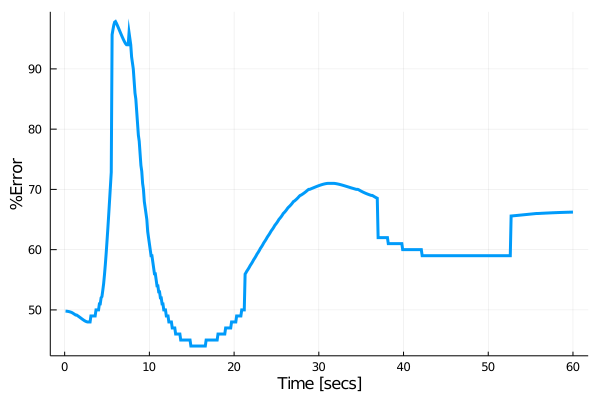}}
\caption{Percentage error history during the reconfiguration manuever.}
\label{f:error_traj}
\end{figure}

With \eqref{e:command_hist}, the trajectory of the relative positions is shown in Fig. \ref{f:traj}.
The force command history from \eqref{e:command_hist} is shown in Fig. \ref{f:f_comtraj}.
The charges computed using Algorithm \ref{algo:bro} are shown in Fig. \ref{f:charge_traj}. 
The minimum percentage error (defined in \eqref{e:perc}) from Algorithm \ref{algo:bro} at each instance of time is plotted in Fig. \ref{f:error_traj}. The average percentage error over the 60 secs reconfiguration maneuver is 63.4\%. This means that the average propellant reduction for the reconfiguration maneuver is 38.6\% compared to using solely thrusters for reconfiguration. 
\section{Conclusions}
In this paper, the control allocation problem for hybrid Coulomb spacecraft formations was defined. The control allocator tries to minimize propellant use by maximizing the contribution from Coulomb forces. The proposed control allocation algorithm uses a minimum matrix rank optimization problem which can be solved efficiently using a trace heuristic. Two numerical examples were provided to demonstrate the proposed control allocator. The first explained the proposed control allocation algorithm and the other demonstrated the propellant savings for a three-spacecraft reconfiguration maneuver in deep space. 

A suggested direction for future work is the design of path planning algorithms that optimize the sequence of relative force commands so as to minimize propellant with a control allocator in the loop. Further work should also be dedicated to improving the algorithm for real-time implementation by using techniques such as warm-starting.  
\bibliography{Coulomb}

\end{document}